\documentclass[12pt]{article}

 \usepackage{a4}
 \usepackage{amsmath}
 \usepackage{amssymb}
\usepackage{graphicx}

\title{Derivation of an integral of Boros and Moll via  convolution of Student t-densities}
\author{Christian Berg  and Christophe Vignat}
 \begin{document}
\maketitle
\begin{abstract} We show that the evaluation of an integral considered
  by Boros and Moll is a special case of a convolution result about
  Student t-densities obtained  by the authors in 2008.
\end{abstract}

\textbf{Keywords} Quartic integral, Student t-density.

\textbf{AMS  Classification Numbers} 33C05, 60E07

\section{Introduction}

In a series of papers \cite{B:M},\cite{A:M},\cite{M:M},\cite{M},\cite{A:M:V} Moll
and his coauthors have considered the integral
\begin{equation}\label{eq:Moll0}
\int_0^\infty\frac{dx}{(x^4+2ax^2+1)^{m+1}}, a>-1,m=0,1,\ldots.
\end{equation}
It was evaluated first by George Boros, who gave the identity
\begin{equation}\label{eq:Moll1}
\int_0^\infty\frac{dx}{(x^4+2ax^2+1)^{m+1}}=\frac{\pi}2\frac{P_m(a)}{[2(a+1)]^{m+1/2}},
\end{equation}
where
\begin{equation}\label{eq:Moll2}
P_m(a)=\sum_{j=0}^m d_{j,m}a^j
\end{equation}
and
\begin{equation}\label{eq:Moll3}
d_{j,m}=2^{-2m}\sum_{i=j}^m 2^{i}\binom{2m-2i}{m-i}\binom{m+i}{m}\binom{i}{j}.
\end{equation}

The paper \cite{A:M} gives a survey of  different proofs of the
formula \eqref{eq:Moll1}.

The purpose of the present paper is to point out that the evaluation
can be considered as a special case of a convolution result about
Student t-densities, thereby adding yet another proof to the list of
\cite{A:M}.

For $\nu>0$  the probability density on $\mathbb R$
\begin{equation}\label{eq:student}
f_{\nu}(x)=\frac{A_{\nu}}{(1+x^2)^{\nu+\tfrac12}},\quad A_{\nu}=
\frac{\Gamma(\nu+\tfrac12)}{\Gamma(\tfrac12)\Gamma(\nu)}
\end{equation}
is called  a Student t-density with $f=2\nu$ degrees of freedom.

It is the special case $\nu=m+1/2$ which is relevant in connection
with the integral \eqref{eq:Moll0}.

The relevant convolution result from \cite{B:V} is 

\begin{equation}\label{eq:conv}
\frac{1}{a}f_{n+\tfrac12}\left(\frac{x}{a}\right)*\frac{1}{1-a}
f_{m+\tfrac12}\left(\frac{x}{1-a}\right)=\sum_{k=n\wedge m}^{n+m}\beta_k^{(n,m)}(a)f_{k+\tfrac12}(x),
\end{equation}
where $0<a<1$, $n,m$ are nonnegative integers and $*$ is the ordinary convolution of densities.

The important issue in \cite{B:V} is to prove that the coefficients
$\beta_k^{(n,m)}(a)$ are non-negative for $0<a<1$. This follows from
explicit formulas for these coefficients in two cases: (I): $n=m$, (II): $n$
arbitrary, $m=0$, combined with the symmetry relation 

\begin{equation}
\beta_{k}^{\left(n,m\right)}\left(a\right)=\beta_k^{\left(m,n\right)}
\left(1-a\right)\label{eq:symmetry}
\end{equation}
and a recursion formula

\begin{equation}
\frac{1}{2k+1}\beta_{k+1}^{\left(n,m\right)}\left(a\right)=\frac{a^{2}}{2n-1}\beta_{k}^{\left(n-1,m\right)}\left(a\right)+\frac{\left(1-a\right)^{2}}{2m-1}\beta_{k}^{\left(n,m-1\right)}\left(a\right).\label{eq:lemma}
\end{equation}

We do not know an explicit formula for  $\beta_k^{(n,m)}(a)$ when
$n,m$ are arbitrary. 
The formula when $m=n$ is given in \cite[Theorem 2.2]{B:V} and reads
\begin{eqnarray}\label{eq:exp1}
\beta_{m+i}^{(m,m)}(a)&=&(4a(1-a))^{i}\left(\frac{m!}{(2m)!}\right)^{2}2^{-2m}
\frac{(2m-2i)!(2m+2i)!}{(m-i)!(m+i)!}\\
&\times&\sum_{j=0}^{m-i}
\binom{2m+1}{2j}\binom{m-j}{i}(2a-1)^{2j},\quad i=0,\ldots,m.
\end{eqnarray}
The case $a=1/2$ leads to
\begin{equation}\label{eq:exp2}
\beta_{m+i}^{(m,m)}(1/2)=\left(\frac{m!}{(2m)!}\right)^{2}2^{-2m}
\frac{(2m-2i)!(2m+2i)!}{(m-i)!(m+i)!}\binom{m}{i}.
\end{equation}

Let us consider $n=m$ and $a=1/2$ in \eqref{eq:conv}, where we replace
$x$ by $x/2$ and multiply by $1/2$ on both sides:
\begin{equation}\label{eq:help1}
f_{m+1/2}\ast f_{m+1/2}(x)=\sum_{k=m}^{2m}\frac12\beta_k^{(m,m)}(1/2)
f_{k+1/2}(x/2).
\end{equation}
The left-hand side is equal to
$$
L:=A_{m+1/2}^2\int_{-\infty}^{\infty}\frac{dy}{[(1+y^2)(1+(x-y)^2)]^{m+1}}
$$
$$
=A_{m+1/2}^2\int_{-\infty}^{\infty}\frac{dt}{[(1+(t+x/2)^2)(1+(t-x/2)^2)]^{m+1}},
$$
where we have used the substitution $t=y-x/2$. Clearly
$$
L=A_{m+1/2}^2\left(1+{x^2}/4\right)^{-2(m+1)}\int_{-\infty}^\infty
\frac{dt}{\left[1+2t^2\frac{1-x^2/4}{(1+x^2/4)^2}+\frac{t^4}{(1+x^2/4)^2}\right]^{m+1}}.
$$
 Finally, substituting
$t=\sqrt{1+x^2/4}\,s$ we get
$$
L=2A_{m+1/2}^2\left(1+{x^2}/4\right)^{-2m-3/2}\int_0^\infty\frac{ds}{\left[1+2as^2+s^4\right]^{m+1}},
$$
where $a=(1-x^2/4)/(1+x^2/4)$.

The right-hand side of \eqref{eq:help1} is equal to
\begin{equation}\label{eq:help2}
R:=\sum_{i=0}^{m}\frac12\beta_{m+i}^{(m,m)}(1/2)\frac{A_{m+i+1/2}}{\left(1+x^2/4\right)^{m+i+1}}.
\end{equation}
Combining this gives
\begin{eqnarray*}
\lefteqn{\int_0^\infty\frac{dx}{(x^4+2ax^2+1)^{m+1}}=}\\
&&\frac{\pi}{2^{2m+2}}\frac{(1/2)_m}{((2m)!)^2}
\sum_{i=0}^m\frac{(2m-2i)!(2m+2i)!}{(m+1/2)_i(m-i)!}\binom{m}{i}\left(1+x^2/4\right)^{m+1/2-i}.
\end{eqnarray*}

Using that $2(a+1)=4/(1+x^2/4)$ we get
$$
\int_0^\infty\frac{dx}{(x^4+2ax^2+1)^{m+1}}=\frac{\pi}{2}\frac{P_m(a)}{[2(a+1)]^{m+1/2}},
$$
where
$$
P_m(a)=\frac{(1/2)_m}{((2m)!)^2}
\sum_{i=0}^m\frac{(2m-2i)!(2m+2i)!}{(m+1/2)_i(m-i)!}\binom{m}{i}\left[(a+1)/2\right])^{i}.
$$
Using the binomial formula for $(a+1)^{i}$ and interchanging the
summations, we finally get
$$
P_m(a)=\sum_{j=0}^m d_{j,m}a^{j}
$$
with 
$$
d_{j,m}=\frac{(1/2)_m}{((2m)!)^2}
\sum_{i=j}^m\frac{(2m-2i)!(2m+2i)!}{(m+1/2)_i(m-i)!2^{i}}\binom{m}{i}\binom{i}{j},
$$
which can easily be reduced to \eqref{eq:Moll3}.

Christian Berg, Department of Mathematical Sciences, University of
Copenhagen, Universitetsparken 5, DK-2100, Copenhagen, Denmark\\
email: berg@math.ku.dk

Christophe Vignat, Laboratoire des Signaux et Syst{\`e}mes, Universit{\'e} d'Orsay, France\\
email: christophe.vignat@u-psud.fr

\begin{thebibliography}{10}

\bibitem{A:M} T.~Amdeberhan and V.~H.~Moll, A formula for a quartic
  integral: a survey of old proofs and some new ones. {\it Ramanujan
    J.} {\bf 18} (2009), 91--102.

\bibitem{A:M:V}  T.~Amdeberhan, V.~H.~Moll and C.~Vignat,  The
  Evaluation of a quartic Integral via Schwinger, Schur and Bessel. Manuscript.
ArXiv:1009.2399v1[math.CA]

\bibitem{B:V} C.~Berg and  C.~Vignat, Linearization coefficients of
  Bessel polynomials and properties of Student
  t-distributions. {\it Const. Approx.} {\bf 27} (2008), 15--32. 

\bibitem{B:M} G.~Boros and V.~H.~Moll, An integral hidden in
  Gradshteyn and Ryzhik. {\it J. Comput. Appl. Math.} {\bf 106}
  (1999), 361--368.

\bibitem{M:M} D.~V.~Manna and V.~H.~Moll, A remarkable sequence of
  integers. {\it Expo Math.} {\bf 27} (2009), 289--312.

\bibitem{M} V.~H~Moll, Seized opportunities, {\it Notices Amer. Math.
    Soc.} {\bf 57} (2010), 476--484.
\end{thebibliography}
\end{document}